\documentclass[a4paper]{amsart}  

\usepackage{amssymb} \usepackage{amscd} \usepackage{times}

\def\zbb{\mathbb{Z}}  
  
  \def\phi{\varphi}
 \def\p1{{\mathbb{P}^1_\zbb}}

\begin{document}

\title{ Note on the Chen-Lin result with Li-Zhang method. }
\author{Samy Skander Bahoura}

\address{}

\email{samybahoura@yahoo.fr}

\date{}

\maketitle

{\bf Abstract:} We give a new proof of Chen-Lin result with Li-Zhang method.

\bigskip

\begin{center}1. INTRODUCTION AND RESULTS.
\end{center}

\bigskip

We set $ \Delta =-\partial_{11} - \partial_{22} $ the geometric Laplacian on $ {\mathbb R}^2 $.

\bigskip

On an open set $ \Omega  $ of $ {\mathbb R}^2 $, with a smooth boundary, we consider the following problem:

$$ (P) \,\,\,\begin{cases}

     \Delta u = V e^u \,\,\, \text{in} \,\,\, \Omega. \\
     
     0 < a \leq  V \leq b <+\infty  \\
               
        \end{cases} $$

\bigskip

The previous equation is called, the Prescribed Scalar Curvature, in
relation with conformal change of metrics. The function $ V $ is the
prescribed curvature.

\bigskip

Here, we try to find some a priori estimates for sequences of the
previous problem.

\smallskip

Equations of the previous type were studied by many authors. We can
see in [B-M], different results for the solutions of those type of
equations with or without boundaries conditions and, with minimal
conditions on $ V $, for example we suppose $ V \geq 0 $ and  $ V
\in L^p(\Omega) $ or $ Ve^u \in L^p(\Omega) $ with $ p \in [1,
+\infty] $. 

We can see in [B-M] the following important Theorem,

\smallskip

{\bf Theorem A}{\it (Brezis-Merle)}.{\it If $ (u_i)_i $ and $ (V_i)_i $ are two sequences of functions relatively to the problem $ (P) $ with, $ 0 < a \leq V_i \leq b < + \infty $, then, for all compact set $ K $ of $ \Omega $,

$$ \sup_K u_i \leq c = c(a, b, m, K, \Omega) \,\,\, {\rm if } \,\,\, \inf_{\Omega} u_i \geq m. $$}

A simple consequence of this theorem is that, if we assume $ u_i = 0 $ on $ \partial \Omega $ then, the sequence $ (u_i)_i $ is locally uniformly bounded.

\bigskip

If, we assume $ V $ with more regularity, we can have another type of estimates, $ \sup + \inf $. It was proved, by Shafrir, see [S], that, if $ (u_i)_i, (V_i)_i $ are two sequences of functions solutions of the previous equation without assumption on the boundary and, $ 0 < a \leq V_i \leq b < + \infty $, then we have the following interior estimate:

$$ C(a/b) \sup_K u_i + \inf_{\Omega} u_i \leq c=c(a, b, K, \Omega). $$

We can see in [C-L], an explicit value of $ C(a/b) =\sqrt {\dfrac{a}{b}} $.

\bigskip

Now, if we suppose $ (V_i)_i $ uniformly Lipschitzian with $ A $ the
Lipschitz constant, then, $ C(a/b)=1 $ and $ c=c(a, b, A, K, \Omega)
$, see Br\'ezis-Li-Shafrir [B-L-S]. This result was extended for
H\"olderian sequences $ (V_i)_i $ by Chen-Lin, see  [C-L]. Also, we
can see in [L], an extension of the Brezis-Li-Shafrir to compact
Riemann surface without boundary. We can see in [L-S] explicit form,
($ 8 \pi m, m\in {\mathbb N}^* $ exactly), for the numbers in front of
the Dirac masses, when the solutions blow-up.

\bigskip

On open set $ \Omega $ of $ {\mathbb R}^2 $ we consider the following equation:

$$ \Delta u_i = V_i e^{u_i} \,\,\, {\rm on } \,\, \Omega. $$

$$ 0 < a \leq V_i \leq b < + \infty, \,\,\, |V_i(x)-V_i(y)|\leq A |x-y|^s, \,\,\, 0 < s < 1, \,\, x, y\in \Omega. $$

{\bf Theorem B}({\it (Chen-Lin}). For all compact $ K \subset \Omega $ and all $ s \in ]0,1[ $ there is a constant $ c= c(a, b, A, s, K, \Omega) $ such that,

$$ \sup_K u_i + \inf_{\Omega} u_i \leq c, \,\, \forall \,\, i. $$

Here we try to prove the previous theorem by the moving-plane method and Li-Zhang method.

\bigskip

We argue by contradiction, and we want to proof that:

$$ \exists \,\, R >0, \,\,\, {\rm such \,\, that} \,\, 4 \log R + \sup_{B_R(0)} u + \inf_{B_{2R}(0)} u \leq c=c(a, b, A), $$

Thus, by contardition we can assume:

$$ \exists \,\, (R_i)_i, \,\, (u_i)_i\,\, R_i \to 0, \,\, 4 \log R_i + \sup_{B_{R_i}(0)} u_i + \inf_{B_{2R_i}(0)} u_i \to +\infty,  $$

\bigskip

\underbar {\it The blow-up analysis}

\bigskip

Let $ x_0 \in \Omega $, we want to prove the theorem locally around $ x_0 $, we use the previous assertion with $ x_0=0 $. The classical blow-up analysis gives the existence of the sequence $ (x_i)_i $ and a sequence of functions $ (v_i)_i $ such that:

\bigskip

We set,

$$ \sup_{B_{R_i(0)}} u_i = u_i(\bar x_i), $$

$$ s_i(x)= 2 \log (R_i-|x-\bar x_i|) + u_i(x), \,\, {\rm and} \,\, s_i(x_i)= \sup_{B_{R_i(\bar x_i)}} s_i, \,\,\, l_i= \dfrac{1}{2}(R_i-|x_i-\bar x_i|). $$

Also, we set:

$$ v_i(x)= u_i[x_i+xe^{-u_i(x_i)/2}]-u_i(x_i), \,\,\, \bar V_i(x) = V_i[x_i+xe^{-u_i(x_i)/2}], $$

Then,

$$ \Delta v_i = \bar V_i e^{v_i}, $$

$$ v_i \leq 2\log 2, \,\, v_i(0)=0, $$

$$ v_i \to v = \log \dfrac{1}{(1+[V(0)/8] |x|^2)^2}\,\, {\rm converge} \,\, {\rm uniformly \,\, on \,\, each \,\, compact \,\, set \,\, of \,\, {\mathbb R}^2 } $$

with $ V(0) = \lim_{i \to + \infty} V_i(x_i) $.

\bigskip

The classical elliptic estimates and the classical Harnack inequality, we can prove the previous uniform convergence on each compact of $ {\mathbb R}^2 $.

\bigskip

\underbar {\it The Kelvin transform and the moving-plane method: Li-Zhang method.}

\bigskip

For $ 0 < \lambda < \lambda_1 $, we define:

$$ \Sigma_{\lambda} = B(0, l_iM_i)-B(0,\lambda). $$

First, we set :

 $$ \bar v_i^{\lambda} = v_i^{\lambda} - 4 \log |x| + 4 \log \lambda = v_i \left ( \dfrac{\lambda^2x}{|x|^2} \right ) + 4 \log \dfrac{\lambda}{|x|}, $$

$$ \bar V_i^{\lambda} = \bar V_i \left ( \dfrac{ \lambda^2 x}{|x|^2} \right ), $$

$$ M_i = e^{u_i(x_i)/2}, $$

and,

$$ w_{\lambda} = \bar v_i - \bar v_i^{\lambda}. $$

Then,

$$ \Delta \bar v_i^{\lambda} = \bar V_i^{\lambda} e^{\bar v_i^{\lambda}}, $$

$$ \min_{|y| = R_i M_i } \bar v_i^{\lambda}  = u_i(x_i+r\theta)-u_i(x_i)+2u_i(x_i) \geq \inf_{\Omega} u_i + u_i(x_i) \to + \infty, $$

and,

$$ \Delta (v_i-v_i^{\lambda}) = \bar V_i (e^{v_i} - e^{v_i^{\lambda}}) + (\bar V_i - \bar V_i^{\lambda})e^{v_i^{\lambda}}, $$

We have the following estimate:

$$ |\bar V_i - \bar V_i^{\lambda}| \leq A M_i^{-s} |x|^s|1-\dfrac{\lambda^2}{|x|^2}|^s, $$


We take an auxiliary function $ h_{\lambda} $ :

Because, $ v_i(x^{\lambda}) \leq C(\lambda_1) < + \infty $, we have,

$$ h_{\lambda} =C_1  M_i^{-s} {\lambda}^2 (\log (\lambda/|x|))+C_2 M_i^{-s} \lambda^{2+s} [1-(\dfrac{\lambda }{|x|})^{2-s}], \,\,\, |x|> \lambda,  $$

with $ C_1, C_2 = C_1, C_2(s, \lambda_1) >0 $,

\bigskip

$$ h_{\lambda} =M_i^{-s} {\lambda}^2(1-\lambda / |x|)(C_1 \dfrac{\log (\lambda / |x|)}{1-\lambda / |x|}+C'_2),  $$

with, $ C'_2 = C'_2(s, \lambda_1) >0 $. We can choose $ C_1 $ big enough to have $ h_{\lambda} <0 $.

\bigskip

\underbar {\it Lemma 1:} There is an $ \lambda_{k,0} >0 $ small enough, such that, for $ 0 < \lambda < \lambda_{k,0} $, we have:

$$ w_{\lambda} + h_{\lambda} >0. $$

\bigskip

We have,

$$ f(r,\theta)= v_i(r\theta) + 2 \log r, $$

then,

$$ \dfrac{\partial f}{\partial r}(r,\theta) = <\nabla v_i(r\theta)|\theta> + \dfrac{2}{r}, $$

According to the blow-up analysis,

$$ \exists \,\, r_0 >0, \,\, C>0, \,\, |\nabla v_i(r\theta)|\theta>| \leq C, \,\, {\rm for } \,\, 0 \leq r < r_0, $$

Then,

$$ \exists \,\, r_0 >0, \,\, C'>0, \,\, \dfrac{\partial f}{\partial r}(r,\theta) > \dfrac{C'}{r}, \,\, 0 < r < r_0, $$

if $ 0 < \lambda < |y| < r_0 $,

$$ w_{\lambda}(y) + h_{\lambda}(y) = v_i(y) - v_i^{\lambda}(y) + h_{\lambda}(y) > C (\log |y|-\log |y^{\lambda}|) + h_{\lambda}(y),$$

by the definition of $ h_{\lambda} $, we have, for $ C, C_0 >0 $ and $ 0 < \lambda \leq  |y| < r_0 $,

$$ w_{\lambda}(y) + h_{\lambda}(y) > (|y|-\lambda ) [ C \dfrac{\log |y| -\log |y^{\lambda}|}{|y|-\lambda} - \lambda^{1+s}C_0 M_i^s ], $$

but,

$$  |y| - |y^{\lambda}| > |y| -\lambda >0, \,\, {\rm and} \,\, |y^{\lambda}|=\dfrac{\lambda^2}{|y|}, $$

thus,

In the first step of the lemma 2, we have,

$$ v_i \geq \min v_i = C_i^1, \,\, v_i^{\lambda}(y) \leq C_1(\lambda_1, r_0), \,\, {\rm if } \,\, r_0 \leq |y| \leq R_iM_i, $$

Thus, in $ r_0 \leq |y| \leq R_iM_i $ and $ \lambda \leq \lambda_1 $, we have,

$$  w_{\lambda} + h_{\lambda} \geq C_i -4\log \lambda + 4 \log r_0 -C' {\lambda_1}^{2+s} $$

then, if $ \lambda \to 0 $, $ -\log \lambda \to + \infty $, and,

$$ w_{\lambda} +h_{\lambda} > 0, \,\,\, {\rm if} \,\, \lambda < \lambda_0^i, \lambda_0^i \,\, {\rm (small) }, \,\, {\rm and} \,\, r_0 < |y| \leq R_i M_i, $$





\bigskip

By the maximum princple and the Hopf boundary lemma, we have:

\bigskip

\underbar {\it Lemma 2:} Let $ \tilde \lambda_k $ be a positive number such that:

$$ \tilde \lambda_k = \sup \{ \lambda < \lambda_1,\,\,\, w_{\lambda}+ h_{\lambda} >0 \,\,\, {\rm in} \,\,\, \Sigma_{\lambda}. \}. $$

Then,

$$ \tilde \lambda_k = \lambda_1. $$

The blow-up analysis gives the follwing inequality for the boundary condition,

\bigskip

For $ y=|y|\theta =R_iM_i\theta $, we have,

$$ w_{\lambda^i}(|y|=R_iM_i) + h_{\lambda^i}(|y|=R_iM_i) = $$

$$ = u_i(x_i+R_i\theta) -u_i(x_i) -v_i(y^{\lambda^i}, |y|=R_iM_i ) - 4 \log \lambda + 4 \log (R_i M_i) + C(s, \lambda_1)M_i^{-s}\lambda^{2+s}[1-(\dfrac{\lambda}{R_iM_i})^{2-s}],  $$





$$ 4 \log \bar R_i + u_i(x_i) + \inf_{\Omega} u_i \to + \infty, $$




\bigskip

which we can write,

$$ w_{\lambda^i}(|y|=R_iM_i) + h_{\lambda^i}(|y|=R_iM_i) \geq \min_{\Omega} u_i + u_i(x_i) + 4 \log R_i -C(s,\lambda_1) \to + \infty, $$

because, $ 0 < \lambda \leq \lambda_1 $,








\bigskip







\newpage

\bigskip

\underbar{\bf References.}

\bigskip

[B-G] L. Boccardo, T. Gallouet. Nonlinear elliptic and parabolic
equations involving measure data. J. Funct. Anal. 87 (1), (1989),
149-169.

\smallskip

[B-L-S] H. Brezis, YY. Li and I. Shafrir. A sup+inf inequality for
some nonlinear elliptic equations involving exponential
nonlinearities. J.Funct.Anal.115 (1993) 344-358.

\smallskip

[B-M] H. Brezis, F. Merle. Uniform estimates and Blow-up Behavior for Solutions of $ -\Delta u=V(x) e^u $ in two dimension. Commun. in Partial Differential Equations, 16 (8 and 9), 1223-1253(1991).

\smallskip

[C-Li] W. Chen, C. Li. A priori Estimates for solutions to Nonlinear Elliptic Equations. Arch. Rational. Mech. Anal. 122 (1993) 145-157.

\smallskip

[C-L] C-C. Chen, C-S. Lin. A sharp sup+inf inequality for a
nonlinear elliptic equation in ${\mathbb R}^2$. Commun. Anal. Geom.
6, No.1, 1-19 (1998).

\smallskip

[L] YY. Li, Harnack type Inequality, the Methode of Moving Planes.
Commun. Math. Phys. 200 421-444.(1999).

\smallskip

[L-S] YY. Li, I. Shafrir. Blow-up Analysis for Solutions of $
-\Delta u = V e^u $ in Dimension Two. Indiana. Math. J. Vol 3, no 4.
(1994). 1255-1270.

\smallskip

[M-W] L. Ma, J-C. Wei. Convergence for a Liouville equation. Comment. Math. Helv. 76 (2001) 506-514.

\smallskip

[R] W. Rudin. Real and Complex Analysis.

\smallskip

[S] I. Shafrir. A sup+inf inequality for the equation $ -\Delta
u=Ve^u $. C. R. Acad.Sci. Paris S\'er. I Math. 315 (1992), no. 2,
159-164.

\end{document}